\documentclass[12pt]{amsart} %
\usepackage{epsfig}
\usepackage{tikz}
\usepackage{array,amsmath, enumerate,  url, psfrag}
\usepackage{amssymb, amsaddr, fullpage}
\usepackage{graphicx,subfigure}
\usepackage{color}

\newcommand\M{\mathcal{M}}
\newcommand\Dv{Dvo\v{r}\'{a}k}

\newtheorem{theorem}{Theorem}[section]
\newtheorem{lemma}[theorem]{Lemma}

\newtheorem{thm}{Theorem}[section]

\newtheorem{definition}[thm]{Definition}

\title{DP-3-coloring of planar graphs without $4,9$-cycles and two cycles from $\{5,6,7,8\}$}

\author{Runrun Liu$^{1}$, Sarah Loeb$^{2}$,  Martin Rolek$^{2}$,   Yuxue Yin$^{1}$,   Gexin Yu$^{1,2}$}

\address{
$^{1}$\small Department of Mathematics, Central China Normal University, Wuhan, Hubei, China.\\
$^2$\small Department of Mathematics, The College of William and Mary, Williamsburg, VA, 23185, USA.
}

\thanks{The research of the last author was supported in part by the Natural Science Foundation of China (11728102) and the NSA grant H98230-16-1-0316.}

\email{gyu@wm.edu}

\begin{document}
\maketitle

\begin{abstract}
A generalization of list-coloring, now known as DP-coloring, was recently introduced by \Dv{} and Postle~\cite{DP18}.
Essentially, DP-coloring assigns an arbitrary matching between lists of colors at adjacent vertices, as opposed to only matching identical colors as is done for list-coloring.
Several results on list-coloring of planar graphs have since been extended to the setting of DP-coloring~\cite{LL18+,LLYY18+,KO18+,KY18+,SN18+,YY18+}.
We note that list-coloring results do not always extend to DP-coloring results, as shown in~\cite{BK18+}.
Our main result in this paper is to prove that every planar graph without cycles of length $\{4, a, b, 9\}$ for $a, b \in \{6, 7, 8\}$ is DP-$3$-colorable, extending three existing results~\cite{SW07,WS11,WLC08} on $3$-choosability of planar graphs.
\end{abstract}

\section{Introduction}

Graphs in this paper are simple and undirected.
We use $V(G)$, $E(G)$, and $F(G)$, respectively, to represent the vertices, edges, and faces of a graph $G$.
A \emph{coloring} of a graph $G$ is a function $c$ that assigns an element $c(v)$ to each vertex $v \in V(G)$.
A \emph{proper coloring} is a coloring such that $c(u) \neq c(v)$ whenever $uv \in E(G)$.

Vizing~\cite{V76}, and independently Erd\H{o}s, Rubin, and Taylor~\cite{ERT76} introduced list coloring, a generalization of proper coloring.
A \emph{list assignment} $L$ gives each vertex $v$ a set of available colors $L(v)$.
A graph is \emph{$L$-colorable} if it has a proper coloring $c$ with $c(v) \in L(v)$ for every vertex $v$.
A graph is \emph{$k$-choosable} (or \emph{$k$-list-colorable}) if it is $L$-colorable whenever $|L(v)| \ge k$ for each $v \in V(G)$.
The \emph{choosability} (or \emph{list chromatic number}) $\chi_\ell(G)$ of a graph $G$ is the least $k$ such that $G$ is $k$-choosable; the analogue for coloring is the \emph{chromatic number} $\chi(G)$.
In the case that $L(V) = [k]$ for each $v \in V(G)$, any $L$-coloring of $G$ is also a \emph{proper $k$-coloring}, where $[k]$ denotes the set of integers $\{1, 2, \dots, k\}$.
Thus we always have $\chi_\ell(G) \ge \chi(G)$.

More recently, \Dv{} and Postle~\cite{DP18} introduced the following idea of correspondence coloring, which has since become known as DP-coloring. This notion generalizes choosability.

\begin{definition} Let $G$ be a simple graph with $n$ vertices and let $L$ be a list assignment for $G$.
For each $v \in V(G)$, let $L_v = \{v\} \times L(v)$.
For each edge $uv \in E(G)$, let $M_{uv}$ be a matching (possibly empty) between the sets $L_u$ and $L_v$, and let $\M_L = \{M_{uv} : uv \in E(G)\}$, called the \emph{matching assignment}.
Let $G_L$ be a graph that satisfies the following conditions:
\begin{itemize}
\item $V\left(G_L\right) = \cup_{v \in V(G)} L_v$,
\item for each $v \in V(G)$, the set $L_v$ is a clique in $G_L$,
\item if $uv \in E(G)$, then the edges between $L_u$ and $L_v$ are exactly those of $M_{uv}$, and
\item if $uv \notin E(G)$, then there are no edges between $L_u$ and $L_v$.
\end{itemize}
We say that $G$ has an $\M_L$-coloring if $G_L$ contains an independent set of size $n$.
The graph $G$ is \emph{DP-$k$-colorable} if, for any list assignment $L$ with $|L(v)| = k$ for each $v \in V(G)$, the graph is $\M_L$-colorable for every matching assignment $\M_L$.
The least $k$ such that $G$ is DP-$k$-colorable is the \emph{DP-chromatic number} of $G$, denoted $\chi_{DP}(G)$.
\end{definition}

We generally identify the elements of $L_v$ with those of $L(v)$ and refer to the elements as colors.
We will often assume without loss of generality that $L(v) = [k]$ for all $v \in V(G)$, as the existence of an independent set in $G_L$ depends only on the matching assignment $\M_L$.
Suppose $G$ is $\M_L$-colorable and $I$ is an independent set of size $n$ in $G_L$.
Then $|I \cap L_v| = 1$ and we refer to the element $i \in I \cap L_v$ as the color given to $v$.

If $L(v) = [k]$ for each $v \in V(G)$ and $M_{uv} = \{(u, i) (v, i) : i \in [k]\}$ for each $uv \in E(G)$, then an $\M_L$-coloring is exactly a proper $k$-coloring.
Additionally, DP-coloring also generalizes $k$-choosability, even with the restriction that $L(v) = [k]$ for each $v \in V(G)$.
To see this, consider a list assignment $L'$ with $|L'(v)| = k$ for all $v \in V(G)$.
For each vertex $v \in V(G)$, there exists a bijection from the elements of $L'(v)$ to $[k]$, and we simply let $M_{uv}$ be the matching between the colors of $u$ and $v$ that correspond to equal elements of $L'(u)$ and $L'(v)$.
Accounting for relabeling, an $\M_L$-coloring is equivalent to an $L'$-coloring.
Thus, any DP-$k$-colorable graph must be $k$-choosable, and so $\chi_{DP}(G) \ge \chi_\ell(G)$ for all graphs $G$.

One difficulty in the study of list coloring is that some techniques useful in solving coloring problems, such as identifation of vertices, are not feasible in the list coloring setting.
DP-coloring can be used to apply these coloring techniques in some situations.
In the paper introducing DP-coloring, \Dv{} and Postle~\cite{DP18} use identification to prove that planar graphs without cycles of length 4 to 8 are 3-choosable.
However, they impose conditions on the matching assignment $\M_L$, and their proof does not give the analogous result that such graphs are DP-3-colorable.
In their paper, \Dv{} and Postle note that DP-coloring is strictly more difficult than list coloring, in the sense that it is possible for $\chi_{DP}(G) > \chi_\ell(G)$ for some graphs $G$.
In particular, they showed that cycles of even length are 2-choosable, but they are not DP-2-colorable.
In addition, while Alon and Tarsi~\cite{AT92} showed that planar bipartite graphs are 3-choosable, Bernshteyn and Kostochka~\cite{BK18+} provide a bipartite planar graph $G$ with $\chi_{DP}(G) = 4$.

These differences, particularly for even cycles, result in difficulties in extending results from list-coloring to DP-coloring.
However, some proofs for list-coloring do extend to DP-coloring.
For example, \Dv{} and Postle note that Tommassen's proofs~\cite{T94,T95} that $\chi_\ell(G) \le 5$ for planar graphs and $\chi_\ell(G) \le 3$ for planar graphs with no 3-cycles or 4-cycles immediately extend to DP-coloring.

There has been considerable recent interest in extending results for choosability of planar graphs to DP-coloring.
Liu and Li~\cite{LL18+}, Sittitrai and Nakprasit~\cite{SN18+}, Kim and Yu~\cite{KY18+}, and Kim and Ozeki~\cite{KO18+} all extend results on 4-choosability of planar graphs to DP-4-coloring.
Yin and Yu~\cite{YY18+} extend results for 3-choosability to DP-3-coloring, in some cases for a larger class of graphs than the analogous choosability result.
Among other results extending that conditions for 3-choosability hold for DP-3-coloring, Liu, Loeb, Yin, and Yu~\cite{LLYY18+} show that planar graphs with no $\{4, 5, 6, 9\}$-cycles or with no $\{4, 5, 7, 9\}$-cycles are DP-3-colorable.
In this paper, we extend the three previous results for 3-choosability of planar graphs stated in Theorem~\ref{thm1}.

\begin{thm}\label{thm1}
A planar graph $G$ is $3$-choosable if one of the following conditions holds
\begin{itemize}

\item $G$ contains no $\{4,6,7,9\}$-cycles. (Wang, Lu, and Chen~\cite{WLC08})
\item $G$ contains no $\{4,6,8,9\}$-cycles. (Shen and Wang~\cite{SW07})
\item $G$ contains no $\{4,7,8,9\}$-cycles. (Wang and Shen~\cite{WS11})
\end{itemize}
\end{thm}

Our main result is the following.

\begin{thm}\label{awesome}
If $a$ and $b$ are distinct values from $\{6, 7, 8\}$, then every planar graph without cycles of lengths $\{4, a, b, 9\}$ is DP-3-colorable.
\end{thm}


Our proofs use the discharging method, which uses strong induction.  We say a structure is \emph{reducible} if it cannot appear in a minimal counterexample $G$.  The proofs of the results in Theorem~\ref{thm1} rely on the fact that an even cycle with all vertices of degree 3 is reducible.  Such a structure is not necessarily reducible in the setting of DP-coloring.  In Section~\ref{lemma}, we use the lemma about ``near-$(k-1)$-degenerate'' subgraphs from~\cite{LLYY18+} which fills a similar role in our reducible structures.
In Section~\ref{lemma}, we also provide our reducible structures and a lemma about how much charge can be given by large faces in our subsequent discharging arguments.
Section~\ref{sec:threeproofs} provides the proofs for Theorem~\ref{awesome}.
We use different initial charges from the ones in~\cite{LLYY18+}, and provide a new unified set of discharging rules for all three cases. 

\section{Lemmas and a brief discussion of the discharging.}\label{lemma}

Graphs mentioned in this paper are all simple.  A $k$-vertex (resp., $k^+$-vertex, $k^-$-vertex) is a vertex of degree $k$ (resp., at least $k$, at most $k$). The \emph{length} of a face is the number of vertices on its boundary, with repetition included.
A face with length $k$ (resp., at least $k$, at most $k$) is a $k$-face (resp., $k^+$-face, $k^-$-face).
We may also refer to an $(\ell_1, \ell_2, \ldots, \ell_k)$-face, which is a $k$-face $f = v_1 v_2 \ldots v_k$ with facial walk $v_1, v_2, \ldots, v_k$ such that $d(v_i) = \ell_i$.
An $(\ell_1, \ell_2, \ldots, \ell_k)$-path and $(\ell_1, \ell_2)$-edge are defined similarly, and we may replace  $\ell_i$ with $\ell_i^+$ to indicate $d(v_i) \ge \ell_i$.
A $3$-vertex is {\em{triangular}} if it is incident to a $3$-face.

\begin{lemma}\label{minimum}
Let $G$ be a smallest graph (with respect to the number of vertices) that is not DP-$k$-colorable. Then $\delta(G)\ge k$.
\end{lemma}

\begin{proof}
Suppose there is a vertex $v$ with $d(v) < k$.
Any $\mathcal{M}_L$-coloring of $G - v$ can be extended to $G$ since $v$ has at most $d(v)$ elements of $L(v)$ forbidden by the colors selected for the neighbors of $v$, while $|L(v)| = k$.
\end{proof}

Let $H$ be a subgraph of $G$.
For each vertex $v \in H$, let $A(v)$ be the set of vertices (of $G_L$) in $L_v$ that are not matched with vertices in $\cup_{u \in G - H} L_u$.
One may think of $A(v)$ as the colors available at $v$ after coloring $G - H$.

\begin{lemma}\label{near-2-degenerate}\cite{LLYY18+}
Let $k \ge 3$ and $H$ be a subgraph of $G$. If the vertices of $H$ can be ordered as $v_1, v_2, \ldots, v_{\ell}$ such that the following hold
\begin{itemize}
\item[(1)] $v_1v_{\ell}\in E(G)$, and $|A(v_1)|>|A(v_{\ell})|\ge 1$,
\item[(2)] $d(v_{\ell})\le k$ and $v_{\ell}$ has at least one neighbor in $G-H$,
\item[(3)] for each $2\le i\le \ell-1$, $v_i$ has at most $k-1$ neighbors in $G[\{v_1, \ldots, v_{i-1}\}]\cup (G-H)$,
\end{itemize}
then a DP-$k$-coloring of $G-H$ can be extended to a DP-$k$-coloring of $G$.
\end{lemma}


For the remainder of this paper, we will let $G$ denote a minimal counterexample to Theorem~\ref{awesome}.
That is $G$ is a planar graph with no $\{4, a, b, 9\}$-cycles, where $a, b \in \{6, 7, 8\}$ are distinct, such that $G$ is not DP-3-colorable, but any planar graph on fewer than $|V(G)|$ vertices with no $\{4, a, b, 9\}$-cycles is DP-3-colorable.
We now use Lemma~\ref{near-2-degenerate} to provide some reducible configurations we will need in Section~\ref{sec:threeproofs}.

\begin{lemma}\label{lem:reducible}
The graph $G$ does not contain any of the following subgraphs:
\end{lemma}
\begin{figure}[h]
\includegraphics[scale=0.18]{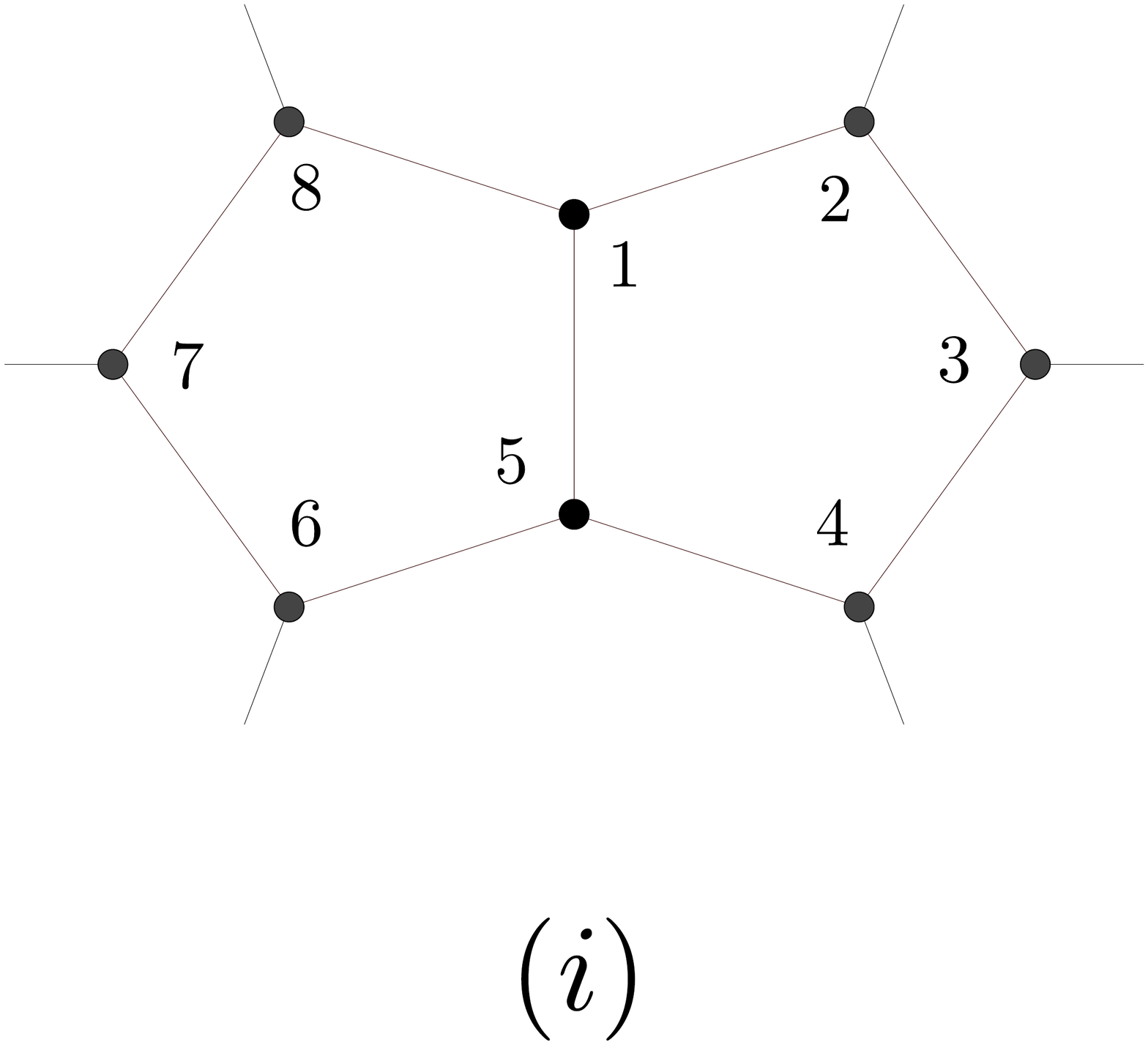} 
\includegraphics[scale=0.17]{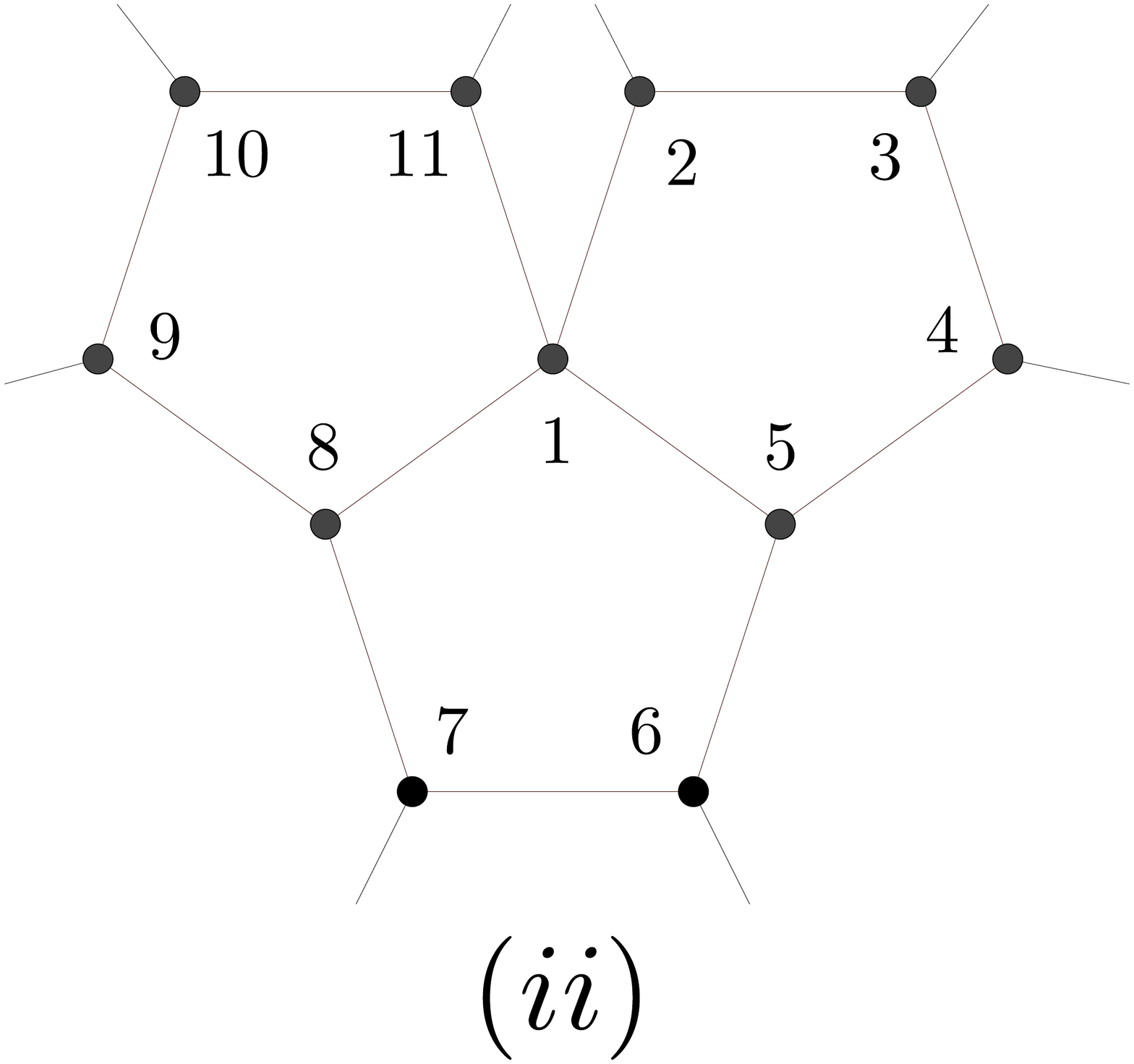}\\
\includegraphics[scale=0.20]{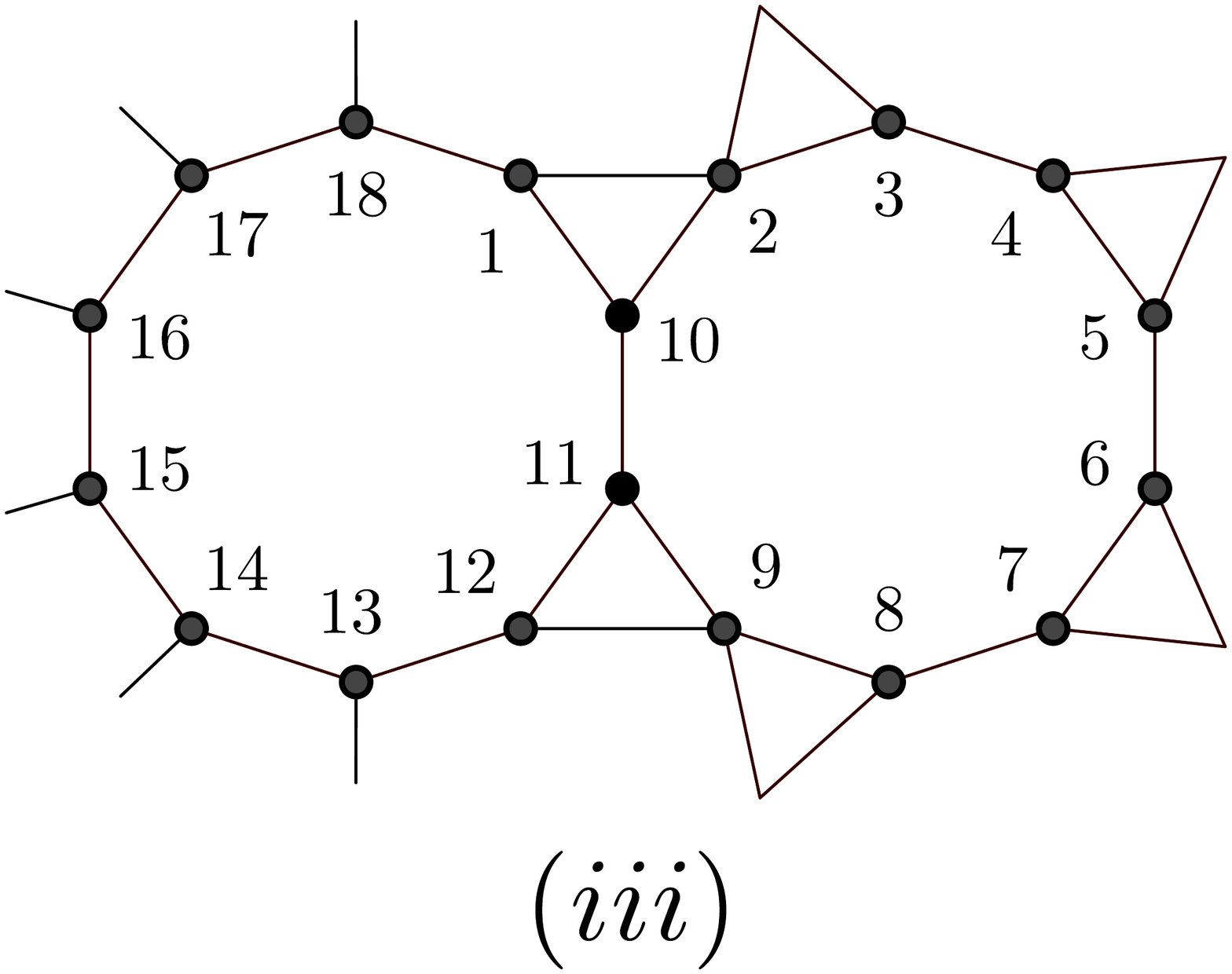}\quad\quad
\includegraphics[scale=0.18]{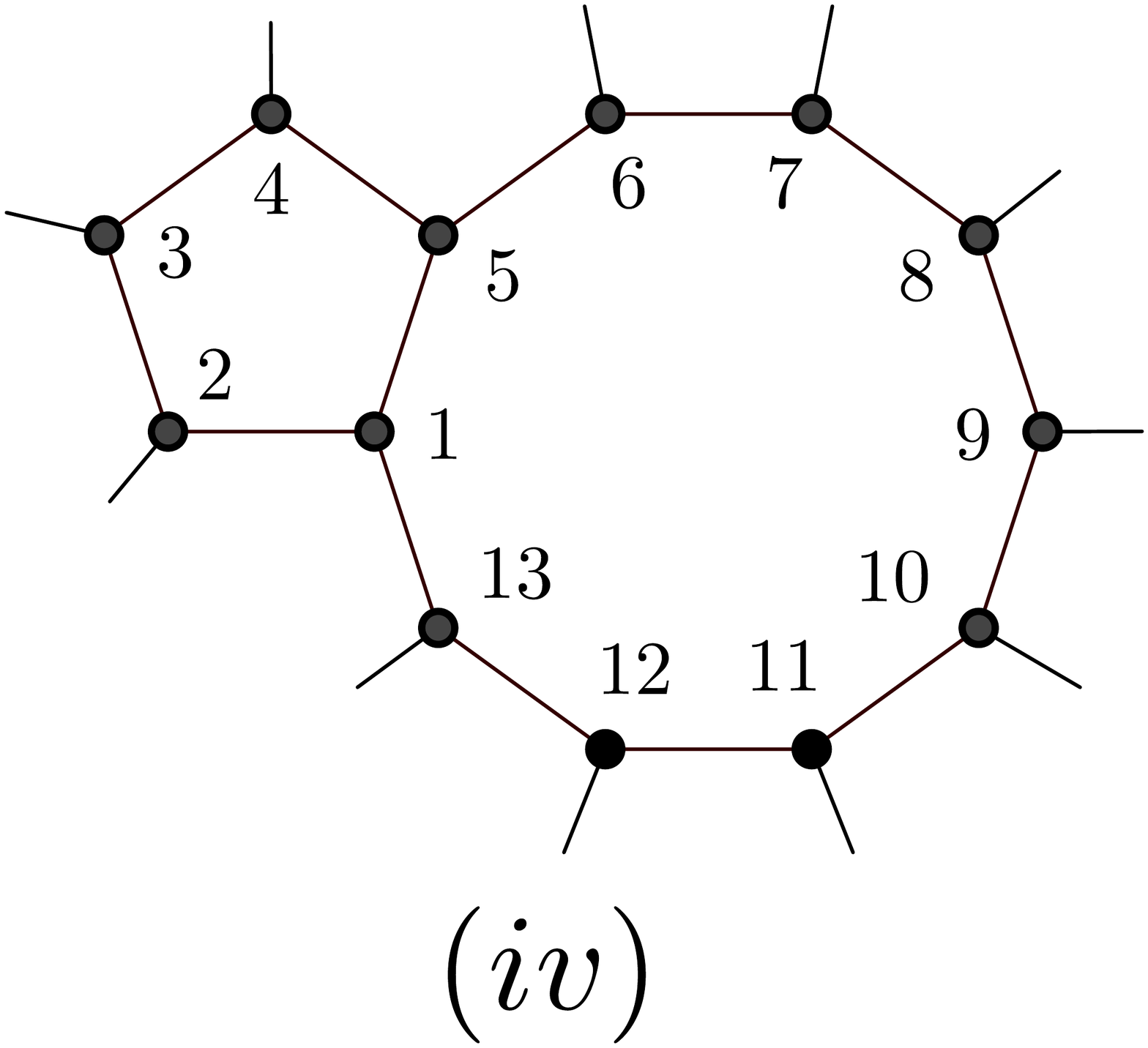}  \quad\quad\includegraphics[scale=0.18]{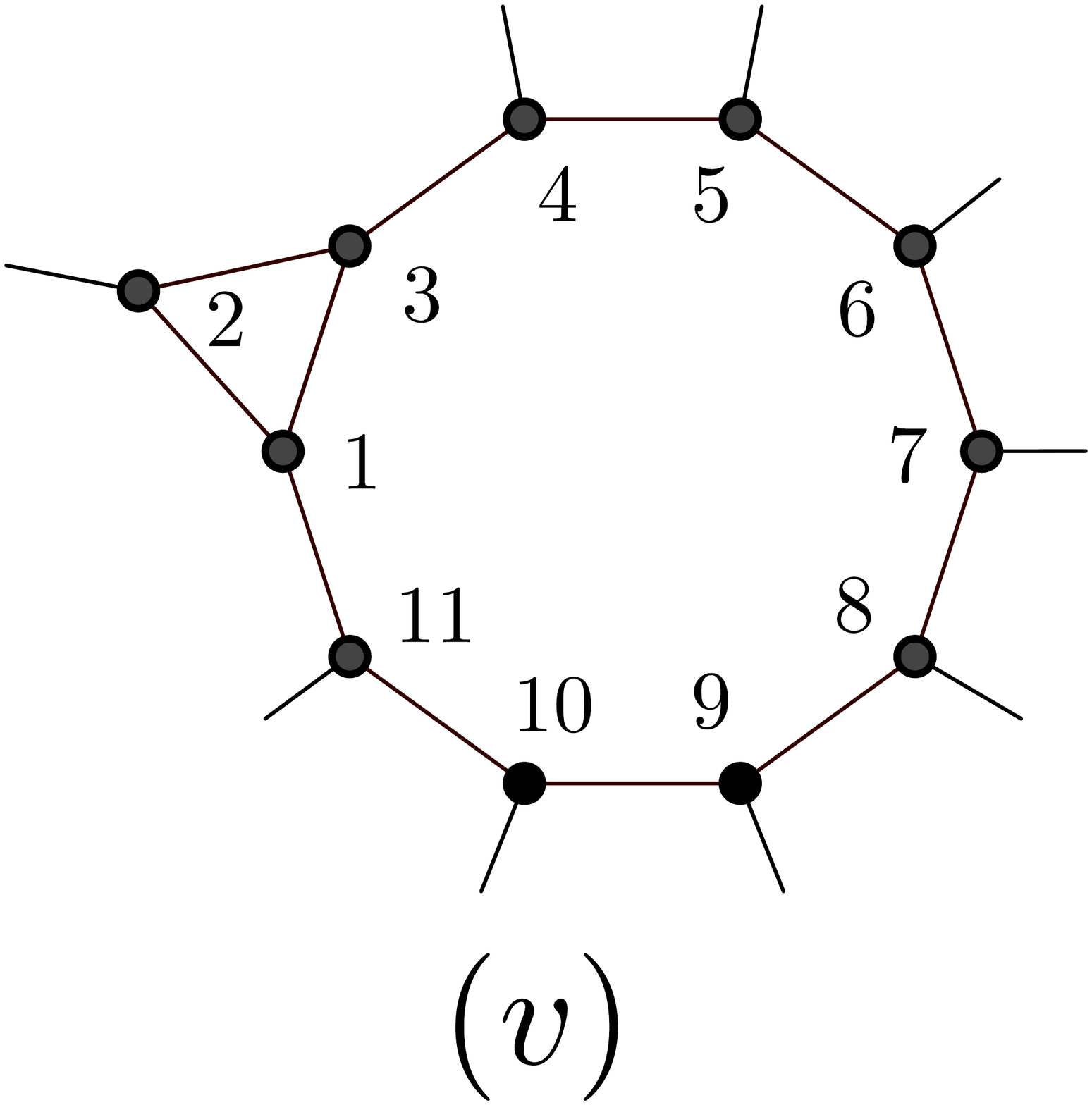}
\end{figure}

\begin{proof}
Let $H$ be the subgraph of $G$ consisting of the labeled vertices, and order the vertices according to their labels.
It is straightforward to verify that all labeled vertices must be distinct, since otherwise cycles of forbidden lengths are created.
From Lemma~\ref{near-2-degenerate}, it follows that a DP-$3$-coloring of $G - H$ can be extended to a DP-$3$-coloring of $G$.
\end{proof}

We use balanced discharging and assign an initial charge of $\mu(x) = d(x) - 4$ to each $x \in V(G) \cup F(G)$.
Let $\mu^*(x)$ be the final charge after the discharging procedure.
From Euler's formula, we have
\begin{equation}\label{eq}
\sum_{x \in V(G) \cup F(G)}(d(x) - 4) = -8.
\end{equation}
We will move charge around and argue in Section~\ref{sec:threeproofs} that each vertex and face ends with non-negative final charge.
This contradiction to (\ref{eq}) will prove our conclusion.


\section{Proof of Theorem~\ref{awesome}} \label{sec:threeproofs}


We say a $10^+$-face $f$ is \emph{good} to a $10$-face $f'$ (and $f'$ is \emph{poor} to $f$) if $f'$ is incident to ten $3$-vertices, and $f$ and $f'$ share a $(3, 3)$-edge such that each end of the shared edge is an end of a $(3, 4^+, 3^+, \dots)$-path of $f$, where the second vertex of the path is incident to either a $(3, 3, 4^+)$-face or a $(3, 3, 3, 3, 4^+)$-face adjacent to both $f$ and $f'$.
A $4^+$-vertex $v$ is {\em{poor}}, {\em{semi-rich}}, or {\em{rich}} to a $10^+$-face $f$ if $v$ is incident to two $3$-faces, one $3$-face, or no $3$-faces adjacent to $f$, respectively.
Moreover, we call a semi-rich $4^+$-vertex $v$ {\em{special}} if $v$ is on a $10^+$-face $f$ such that $v$ is rich to $f$, and we call a $5$-face {\em{bad}} if it is incident to five $3$-vertices and adjacent to two $5$-faces.

The following are our discharging rules:

\begin{enumerate}[(R1)]
\item Each $3$-face gets $\frac{1}{3}$ from each adjacent $5^+$-face, and each $5^+$-face gets $\frac{1}{5}$ from each incident $5^+$-vertex.

\item Each $10^+$-face gets $\frac{1}{6}$ from each incident special semi-rich $4^+$-vertex, and gives $\frac{1}{3}$ to each incident rich $4$-vertex which is on a $3$-face.

\item If $f$ is good to $f'$, then $f$ gives $\frac{1}{6}$ to $f'$.


\item[(R4a)] If $G$ contains no $\{4, 7,8,9\}$-cycles, then each $5$-face that shares a $(3,3)$-edge with a $3$-face gets $\frac{1}{3}$ from each adjacent $10^+$-face, each $5$-face that shares a $(3,4^+)$-edge with a $3$-face gets $\frac{1}{3}$ from the $10^+$-face incident to the $3$-vertex of the $(3,4^+)$-edge, and each $5$-face sends $1$ to any incident triangular $3$-vertex and $\frac{1}{3}$ to any adjacent $3$-face, then distributes its remaining charge evenly to its adjacent $10$-faces. Each $3$-vertex gets its remaining needed charge evenly from its incident $6^+$-faces.    

\item[(R4b)] If $G$ contains no $\{4, 6, a, 9\}$-cycles for $a \in \{7, 8\}$, then each 3-vertex gets $1$ evenly from its incident $5^+$-faces, each $5$-face gets $\frac{1}{6}$ from each adjacent $7^+$-face, and each bad $5$-face gets an additional $\frac{1}{12}$ from each adjacent $5$-face.
Following this, if a $5$-face has positive charge, then  it distributes its surplus charge to its adjacent $5$-faces.
\end{enumerate}

\begin{lemma}
Each vertex and each $8^-$-face have non-negative final charge.
\end{lemma}

\begin{proof}

Note that if a graph $G$ contains no $\{4, 7, 8, 9\}$-cycles, then each $3$-vertex must be incident to at least one $6^+$-face.
By (R4a) and (R4b), the final charge of each $3$-vertex is $0$.
Note that by (R2), a $4$-vertex $v$ sends out charge only if it is special, in which case $v$ is rich to a $10^+$-face and semi-rich to at most two $10^+$-faces.
So if $d(v) = 4$, then $\mu^*(v) \ge (4 - 4) - \frac{1}{6} \cdot 2 + \frac{1}{3} = 0$.
Now let $v$ be a $5^+$-vertex.
Then by (R1) and (R2), $v$ sends out at most $\frac{1}{5}$ to each of its incident faces, so $\mu^*(v) \ge (d(v) - 4) - \frac{1}{5}d(v) \ge 0$.
Hence all vertices end with non-negative charge.

Let $f$ be an $8^-$-face in $G$.
If $d(f) = 3$, then $f$ is adjacent to three $5^+$-faces.
So by (R1), $\mu^*(f) \ge 0$.
For the remaining $8^-$-faces, we consider the consequences of rules (R4a) and (R4b) as separate cases.

{\bf Case 1: $G$ has no $\{4,7,8,9\}$-cycles.}

Then the only $8^-$-faces of $G$ left to consider are $5$- and $6$-faces.
Note that a $5$-face cannot be adjacent to a $6$-face.
If $d(f) = 6$, then $f$ is adjacent to no $3$-face and by (R4a) sends at most $\frac{1}{3}$ to each incident vertex.
So $\mu^*(f) \ge (6 - 4) - \frac{1}{3} \cdot 6 = 0$.
If $d(f) = 5$, then $f$ is adjacent to at most one $3$-face and at least four $10^+$-faces.
If $f$ is not adjacent to any $3$-face, then it only needs to send $1$ evenly to its adjacent 10-faces by (R4a), so $\mu^*(f) \ge 0$.
If $f$ shares a $(3,3)$-edge with a $3$-face, then by (R1) and (R4a), $f$ sends $\frac{1}{3}$ to its adjacent $3$-face and $1$ to each incident triangular $3$-vertex, and $f$ gets $\frac{1}{3}$ from each adjacent $10^+$-face.
So $\mu^*(f) \ge (5 - 4) - \frac{1}{3} - 1 \cdot 2 + \frac{1}{3} \cdot 4 = 0$.
If $f$ shares a $(3,4^+)$-edge with a $3$-face, then by (R1) and (R4a), $f$ sends $\frac{1}{3}$ to its adjacent $3$-face and $1$ to its incident triangular $3$-vertex, and $f$ gets $\frac{1}{3}$ from one $10^+$-face.
So $\mu^*(f) \ge (5 - 4) - \frac{1}{3} - 1 + \frac{1}{3} = 0$.
If $f$ shares a $(4^+, 4^+)$-edge with a $3$-face, then by (R1) and (R4a), $f$ sends $\frac{1}{3}$ to its adjacent $3$-face and sends the remaining charge evenly to adjacent $10$-faces.
Thus, $\mu^*(f)\ge0$.

{\bf Case 2: $G$ has no $\{4, 6, a, 9\}$-cycles for $a \in \{7, 8\}$.}

Let $f$ be an $8$-face.
Then $G$ contains no $\{4, 6, 7, 9\}$-cycles, so $f$ is not adjacent to any 3-face.
Thus, by (R4b), $f$ gives $\frac{1}{3}$ to each incident 3-vertex and $\frac{1}{6}$ to each adjacent 5-face.
So $\mu^*(f) \ge (8 - 4) - \frac{1}{3} \cdot 8 - \frac{1}{6} \cdot 8 = 0$.

Let $f$ be a $7$-face.
Then $G$ contains no $\{4, 6, 8, 9\}$-cycles, so $f$ is again not adjacent to any 3-face.
By (R1) and (R4b), $f$ gives $\frac{1}{3}$ to each incident 3-vertex and $\frac{1}{6}$ to each adjacent 5-face.  Note that no two $5$-faces are adjacent.
So $\mu^*(f) \ge (7 - 4) - \max\big\{\frac{1}{3} \cdot 7 + \frac{1}{6} \cdot 3, \frac{1}{3} \cdot 6 + \frac{1}{6} \cdot 4, \frac{1}{3} \cdot 5 + \frac{1}{6} \cdot 7 \big\} > 0$.

It remains to consider $5$-faces.
So let $f$ be a $5$-face.
Let $r_5(f)$ and $s_3(f)$ be the number of adjacent 5-faces and incident 3-vertices of $f$, respectively.
Since $G$ contains no $6$-cycles, $f$ cannot be adjacent to a $3$-face.

If $f$ is a bad $5$-face, then $s_3(f) = 5$ and $r_5(f) = 2$.
By Lemma~\ref{lem:reducible}$(i)$, $f$ is not adjacent to another bad $5$-face.
Then by (R4b), $f$ gives $\frac{1}{3}$ to each incident $3$-vertex, and $f$ gets $\frac{1}{6}$ from each adjacent $7^+$-face and $\frac{1}{12}$ from each adjacent $5$-face.
So $\mu^*(f) \ge (5 - 4) - \frac{1}{3} \cdot 5 + \frac{1}{12} \cdot 2 + \frac{1}{6} \cdot 3 = 0$.
Thus we may assume that $f$ is not bad.
By (R4b),
\[ \mu^*(f) \ge (5 - 4) - \frac{1}{3} s_3(f) + \frac{1}{6} (5 - r_5(f)) - \frac{1}{12} b_5(f) = \frac{1}{6} \left(11 - 2s_3(f) - r_5(f) - \frac{1}{2}b_5(f)\right), \]
where $b_5(f)$ is the number of bad $5$-faces adjacent to $f$.
Clearly, $b_5(f)\le r_5(f)$.

Note that no $3$-vertex can be incident to three 5-faces since $G$ contains no $9$-cycles. 
If $s_3(f) \le 2$, then $r_5(f) \le 5$ and $b_5(f) \le 1$.
So $\mu^*(f) \ge \frac{1}{6} \left(11 - 2 \cdot 2 - 5 - \frac{1}{2}\right) > 0$.
If $s_3(f) = 3$, then $r_5(f) \le 3$ and $b_5(f) \le 1$.
Thus $\mu^*(f) \ge \frac{1}{6} \left(11 - 2 \cdot 3 - 3 - \frac{1}{2}\right) > 0$.
If $s_3(f) = 5$, then $r_5(f) \le 1$ since $f$ is not bad.
By Lemma~\ref{lem:reducible}$(i)$, $b_5(f) = 0$, so $\mu^*(f) \ge \frac{1}{6}(11 - 2 \cdot 5 - 1) = 0$.
Lastly, if $s_3(f) = 4$, then $r_5(f) \le 3$.
So $\mu^*(f) \ge \frac{1}{6}\left(11 - 2 \cdot 4 - r_5(f) - \frac{1}{2}b_5(f)\right)$, and $\mu^*(f) < 0$ only if $r_5(f) = 3$ and $b_5(f) = 1$, in which case $\mu^*(f) \ge -\frac{1}{12}$.
Let $v$ be the $4^+$-vertex incident to $f$. 
If $d(v) \ge 5$, then by (R1), $f$ gets $\frac{1}{5} > \frac{1}{12}$ from $v$ and ends with non-negative charge.
If $d(v) = 4$, then by Lemma~\ref{lem:reducible}$(ii)$, at least one of the $5$-faces adjacent to $f$ and incident to $v$, say $f'$, has at least two $4^+$-vertices.
Note that $r_5(f') \le 4$.
By (R4b), $f'$ can send
\begin{align*}
&\frac{1}{r_5(f')} \cdot \frac{11 - 2s_3(f') - r_5(f') - \frac{1}{2} b_5(f')}{6} \\
&\ge \max\left\{\frac{1}{4} \cdot \frac{11 - 2 \cdot 2 - 4}{6},\ \frac{1}{3} \cdot \frac{11 - 2 \cdot 3 - 3 - \frac{1}{2}}{6},\ \frac{1}{2} \cdot \frac{11 - 2 \cdot 3 - 2 - \frac{1}{2}}{6},\ \frac{11 - 2 \cdot 3 - 1}{6}\right\}\\
&> \frac{1}{12}
\end{align*}
to each of its adjacent $5$-faces, which includes $f$.
Hence $\mu^*(f) \ge 0$.
\end{proof}

We now only need to verify that $10^+$-faces end with non-negative charge.
Let $f$ be a $10^+$-face.
Let $P$ be a maximal path (or possibly a cycle) along $f$ such that every edge of $P$ is adjacent to a $5^-$-face.
Let $\mathcal{P}$ be a collection of all such paths $P$ along $f$.
By construction, the paths of $\mathcal{P}$ are disjoint.
Let $t_i$ denote the number of paths of $\mathcal{P}$ with $i$ vertices for $i \ge 2$, and let $t_1$ denote the number of vertices incident to $f$ not contained in any path $P \in \mathcal{P}$.
Then $\sum_{i \ge 1} i \cdot t_i = d(f)$.
We will use the following two lemmas to simplify our final analysis.

\begin{lemma}\label{12+faces}
A $10^+$-face $f$ can afford to give out at least \begin{equation}\label{eq-general}
\frac{1}{3} \sum_{i \ge 1} i \cdot t_i + \frac{2}{3} \sum_{i \ge 2} t_i + \frac{1}{3}t_1 + \frac{1}{3} \sum_{i \ge 3} (i - 2)t_i - \frac{x}{3},
\end{equation}
where $x = 0$ if $d(f) \ge 12$, $x = 1$ if $d(f) = 11$, and $x = 2$ if $d(f) = 10$.
\end{lemma}

\begin{proof}
From $\sum_{i \ge 1} i \cdot t_i = t_1 + 2 \sum_{i \ge 2} t_i + \sum_{i \ge 3} (i - 2)t_i$, we have
\[ 2 \sum_{i \ge 2} t_i = \sum_{i \ge 1} i \cdot t_i - t_1 - \sum_{i \ge 3}(i - 2)t_i. \]
Therefore,
\begin{align*}
&\frac{1}{3} \sum_{i \ge 1} i \cdot t_i + \frac{2}{3} \sum_{i \ge 2} t_i + \frac{1}{3} t_1 + \frac{1}{3} \sum_{i \ge 3} (i - 2)t_i - \frac{x}{3} \\
&= \frac{1}{3} \sum_{i \ge 1} i \cdot t_i + \frac{1}{3} \left(\sum_{i \ge 1} i \cdot t_i - t_1 - \sum_{i \ge 3}(i - 2)t_i\right) + \frac{1}{3} t_1 + \frac{1}{3} \sum_{i \ge 3} (i - 2)t_i - \frac{x}{3} \\
&= \frac{2}{3} \sum_{i \ge 1} i \cdot t_i - \frac{x}{3} \\
&= \frac{2}{3}d(f) - \frac{x}{3}.
\end{align*}
So $\mu^*(f) \ge (d(f) - 4) - \frac{2}{3}d(f) + \frac{x}{3} = \frac{1}{3}(d(f) - 12 + x) \ge 0$ when $d(f) \ge 10$.
\end{proof}

\begin{lemma}\label{lem:eq-10+}
Each $10^+$-face $f$ needs to send out at most
\begin{equation}\label{eq-10+}
\frac{1}{3} \sum_{i \ge 1} i \cdot t_i + \frac{2}{3} \sum_{i \ge 2} t_i + \frac{1}{6} \sum_{i \ge 6}(i - 5)t_i + \frac{1}{6} \sum_{i \ge 3}t_i,
\end{equation}
\end{lemma}

\begin{proof}
Note that $f$ gives at most $\frac{1}{3}$ to each vertex not on any path of $\mathcal{P}$ by (R2), (R4a) or (R4b).
Now we show that $f$ needs to send out at most $\sum_{i \ge 2} \frac{1}{3}(i + 2)t_i + \frac{1}{6} \sum_{i \ge 6} (i - 5)t_i$ to the $5^-$-faces and vertices along the paths of $\mathcal{P}$.
In the case of forbidding $\{4, a, 8, 9\}$-cycles with $a \in \{6, 7\}$, $f$ gives at most $\frac{1}{2}$ to each endpoint and $\frac{1}{3}$ to each adjacent $5^-$-face along an $i$-path.
So $f$ gives $\frac{1}{2} \cdot 2 + \frac{1}{3}(i - 1) = \frac{1}{3}(i + 2)$ to each $i$-path.
In the case of forbidding $\{4, 6, 7, 9\}$-cycles, $f$ may instead need to give at most $\frac{1}{3}$ to each vertex and $\frac{1}{6}$ to each adjacent $5$-face along an $i$-path when all vertices of the path are $3$-vertices and all adjacent faces are $5$-faces.
Then $f$ gives $\frac{1}{3} i + \frac{1}{6}(i - 1) = \frac{1}{2} i - \frac{1}{6}$ to each $i$-path of this form.
In addition, $\frac{1}{2} i - \frac{1}{6} > \frac{1}{3}(i + 2)$ only if $i \ge 6$.
So in any case, $f$ needs to send to each $i$-path at most $\sum_{i \ge 2} \frac{1}{3}(i + 2)t_i + \sum_{i \ge 6} \frac{1}{6}(i - 5)t_i$.
By (R3), $f$ may need to send out an additional $\frac{1}{6} \sum_{i \ge 3}t_i$ to poor $10$-faces.
Note that if $f$ sends charge over a $2$-path to a poor face, then the $2$-path has a $4^+$-vertex as an endpoint, and $f$ sends at most $\frac{1}{2} + \frac{1}{3} + \frac{1}{6} \le \frac{1}{3}(2 + 2)$ across this path, so the $\frac{1}{6}$ sent to the poor face is already accounted for in the above formula.
Therefore $f$ sends at most
\[ \frac{1}{3} t_1 + \sum_{i \ge 2} \frac{1}{3}(i + 2)t_i + \sum_{i \ge 6} \frac{1}{6}(i - 5)t_i + \frac{1}{6} \sum_{i \ge 3} t_i \]
to its incident vertices and adjacent faces, from which ~(\ref{eq-10+}) follows.
\end{proof}

Assume that $\mu^*(f) < 0$.
Let $s(f)$ be the number of semi-rich $4$-vertices and $5^+$-vertices on $f$.
Note that each semi-rich $4$-vertex on $f$ saves at least $\frac{1}{3} - \frac{1}{6} = \frac{1}{6}$ and by (R2) each $5^+$-vertex on $f$ gives $\frac{1}{5}$ to $f$.
Then by Lemmas~\ref{12+faces} and~\ref{lem:eq-10+},
\[ \frac{1}{6}\left(2t_1 + 2\sum_{i \ge 3} (i - 2)t_i - 2x\right) < \frac{1}{6} \left(\sum_{i \ge 6} (i - 5)t_i + \sum_{i \ge 3} t_i - s(f)\right), \]
which implies that
\[ 2x > s(f) + 2t_1 + \sum_{i \ge 3} (2i - 5)t_i - \sum_{i \ge 6} (i - 5)t_i \ge s(f) + 2t_1 + t_3 + 3t_4 + 5 \sum_{i \ge 5} t_i. \]

Clearly, $x > 0$, so $d(f) \le 11$.
Recall that $f$ gives at most $\frac{1}{3}(i + 2)$ across any $i$-path with $i \le 5$, and note that this is only possible when the ends of the path are $3$-vertices requiring $\frac{1}{2}$ from $f$, and each face adjacent to $f$ along the path is a $5^-$-face requiring $\frac{1}{3}$ from $f$.
Let $d(f) = 11$.
Then $x = 1$.
So $t_1 = 0$, $t_3 \le 1$ and $t_i = 0$ for $i \ge 4$.
By parity, $t_3 = 1$ and $t_2 = 4$, so $s(f) = 0$.
In this case, $f$ is not good to any $10$-face, so we have $\mu^*(f) \ge (11 - 4) - \frac{5}{3} - 4 \cdot \frac{4}{3} = 0$.

Let $d(f) = 10$.
Then $x = 2$, and $t_1 \le 1$ and $t_i = 0$ for $i \ge 5$.
Let $t_1 = 1$.
Then $t_3 \le 1$ and $t_4 = 0$.
By parity, $t_3 = 1$, and $t_2 = 3$.
It follows that $s(f) = 0$.
In this case, $f$ is not good to any $10$-face, so $\mu^*(f) \ge (10 - 4) - \frac{5}{3} - 3 \cdot \frac{4}{3} - \frac{1}{3} = 0$.
Thus we may assume $t_1 = 0$.
Then $s(f) + t_3 + 3t_4 \le 3$.
By parity, we have three primary cases: $t_4 = 1$ and $t_2 = 3$, or $t_3 = t_2 = 2$, or $t_2 = 5$.

In the first case, $s(f) = 0$, so $f$ is not good to any $10$-face, and $\mu^*(f) \ge 6 - 2 - 3 \cdot \frac{4}{3} = 0$.

In the second case, $s(f) \le 1$, and $f$ is good to at most two $10$-faces.
If $s(f) = 0$, then by Lemma \ref{lem:reducible}$(iii)$, $f$ cannot be good to a $10$-face, so $\mu^*(f) \ge 6 - 2 \cdot \frac{5}{3} - 2 \cdot \frac{4}{3} = 0$.
So let $s(f) = 1$.
If $f$ is good to at most one $10$-face, then $\mu^*(f) \ge 6 - 2 \cdot \frac{5}{3} - 2 \cdot \frac{4}{3} - \frac{1}{6} + \frac{1}{6} = 0$, where the final $\frac{1}{6}$ is the minimum $f$ saves or receives from the vertex counted by $s(f)$.
If $f$ is good to two $10$-faces, then the $4^+$-vertex counted by $s(f)$ must be the end of a $2$-path, so $f$ saves at least $\frac{1}{3}$ from this vertex and $\mu^*(f) \ge 6 - 2 \cdot \frac{5}{3} - 2 \cdot \frac{4}{3} - 2 \cdot \frac{1}{6} + \frac{1}{3} = 0$.

In the last case, $\mu^*(f) \ge 6 - \frac{4}{3} \cdot 5 = -\frac{2}{3}$.
We first assume that $G$ contains no $\{4, 7, 8, 9\}$-cycles.
Note that by (R4a), $f$ gives no charge to adjacent special $5$-faces, where a $5$-face is special if it does not share a $(3,3)$-edge with a $3$-face.
Thus we may assume that $f$ is adjacent to at most one special $5$-face, for otherwise, $f$ saves at least $\frac{2}{3}$ and ends with non-negative charge.
If $f$ is incident to at least two $4^+$-vertices, then $f$ saves at least $2\left(\frac{1}{2} - \frac{1}{6}\right) = \frac{2}{3}$, where the $\frac{1}{6}$ is because $f$ may now be good to adjacent $10$-faces.
If $f$ is incident to one $4^+$-vertex, then $f$ saves at least $\frac{1}{2} + \frac{1}{6} = \frac{2}{3}$ since the $4^+$-vertex must be rich to a $10^+$-face adjacent to $f$, and $f$ cannot be good to any $10$-face.
We may thus assume that $f$ is incident to ten $3$-vertices.
By Lemma~\ref{lem:reducible}$(iv)$ and $(v)$, each $5^-$-face adjacent to $f$ must contain a $4^+$-vertex.
This implies that all $3$-faces adjacent to $f$ are $(3,3,4^+)$-faces, and each adjacent $5$-face other than at most one special $5$-face contains a $4^+$-vertex and shares a $(3,3)$-edge with a $3$-face.
It follows that $f$ is poor to at least three $10^+$-faces if $f$ contains a special $5$-face, and is poor to five $10^+$-faces if $f$ contains no special $5$-faces.
Therefore, by (R2) and (R4a), $f$ receives $\min\left\{\frac{1}{3} + \frac{1}{6} \cdot 3, \frac{1}{6} \cdot 5\right\} > \frac{2}{3}$ from adjacent $10^+$-faces it is poor to, so $f$ ends with non-negative charge.
Now we assume that $G$ contains no $\{4, 6, a, 9\}$-cycles for $a \in \{7, 8\}$.
If $f$ is adjacent to a $5$-face, then $f$ gives at most $\frac{1}{3} \cdot 2 + \frac{1}{6} = \frac{5}{6}$ across this $2$-path.
Thus we may assume that $f$ is adjacent to at most one $5$-face, for otherwise, $f$ saves $2\left(\frac{4}{3} - \frac{5}{6}\right) > \frac{2}{3}$ and ends with non-negative charge.
If $f$ is adjacent to exactly one $5$-face and at least one $4^+$-vertex, then $f$ saves at least $\left(\frac{4}{3} - \frac{5}{6}\right) + \left(\frac{1}{3} - \frac{1}{6}\right) = \frac{2}{3}$.
If $f$ is not adjacent to any $5$-faces, then $f$ saves at least $2\left(\frac{1}{2} - \frac{1}{6}\right) = \frac{2}{3}$ if $f$ is incident to at least two $4^+$-vertices, and at least $\frac{1}{2} + \frac{1}{6} = \frac{2}{3}$ if $f$ is incident to exactly one $4^+$-vertex, where in the latter case the $4^+$-vertex must be a special semi-rich $4$-vertex which gives $\frac{1}{6}$ to $f$ by (R2).
We may therefore assume that all vertices incident to $f$ are $3$-vertices.
By Lemma~\ref{lem:reducible}$(v)$, all $3$-faces adjacent to $f$ are $(3,3,4^+)$-faces.
Thus $f$ is poor to five $10^+$-faces if $f$ is not adjacent to a $5$-face, and $f$ is poor to three $10^+$-faces otherwise.
Therefore $f$ gets $\frac{1}{6}$ from each $10^+$-face good to $f$ and saves $\frac{4}{3} - \frac{5}{6}$ if it is adjacent to a $5$-face, for a total of at least $\min\left\{\frac{1}{6} \cdot 5, \left(\frac{4}{3} - \frac{5}{6}\right) + \frac{1}{6} \cdot 3 \right\} > \frac{2}{3}$.
Therefore in all cases, $f$ ends with non-negative charge.

\section{Final remarks.}

We remark that we are yet unable to prove that planar graphs without $\{4, 5, 8, 9\}$-cycles are DP-3-colorable.  While some of our lemmas may be useful in such a proof (in particular, Lemmas~\ref{near-2-degenerate} and~\ref{12+faces}), this case is considerably more difficult than any of the three cases of Theorem~\ref{awesome}, and a unified proof of all four cases does not seem possible. Another remark is that \Dv{} and Postle showed that planar graphs without cycles of lengths from $4$ to $8$ are ``weakly'' DP-3-colorable. It remains open to know if such planar graphs are DP-$3$-colorable.

\end{document}